%% file: main.tex
\documentclass[prodmode,acmtoms]{acmsmall}
\input{preamble}

\input{macros}

\usepackage{natbib}
\usepackage{setspace}
\begin{document} % Actual document begins here

%\setstretch{2.5}

\input{frontmatter}
\maketitle

\input{introduction}
\input{background}
\input{parallelization-scheme}

\input{structure-of-the-code}
\input{numerical-examples}

\input{conclusion}

% Acknowledgments
\begin{acks}
The computations presented in Sec.~\ref{sec:num} were made possible by the facilities of the Shared Hierarchical 
Academic Research Computing Network (SHARCNET:www.sharcnet.ca) and Compute/Calcul Canada.
\end{acks}

% Bibliography
\bibliographystyle{acmsmall}
\bibliography{pampac}

% History dates
\received{September 2013}{}{}

\end{document}

%% file: preamble.tex
\usepackage{amsmath,amssymb}
\usepackage{xspace}
\usepackage[algoruled,vlined,linesnumbered]{algorithm2e}

\SetAlFnt{\small}
\SetAlCapFnt{\small}
\SetAlCapNameFnt{\small}
\SetAlCapHSkip{0pt}
\IncMargin{-\parindent}
\SetKwComment{Comment}{\%\ }{}
\DontPrintSemicolon

%% file: macros.tex
\renewcommand{\hat}{\widehat}
\renewcommand{\vec}[1]{{\mathbf{#1}}}

\newcommand{\norm}[1]{\ensuremath{\left\lVert #1 \right\rVert}}
\newcommand{\R}{\ensuremath{\mathbb{R}}}

\renewcommand{\Re}{\ensuremath{\operatorname{Re}}\xspace}

\newcommand{\lambi}{\ensuremath{\lambda_{\mathrm{min}}}}
\newcommand{\lambf}{\ensuremath{\lambda_{\mathrm{max}}}}
\newcommand{\ie}{\textit{i.e.}}
\newcommand{\eg}{\textit{e.g.}}
\newcommand{\pampac}{\textsc{PAMPAC}\xspace}
\newcommand{\RED}{\texttt{RED}\xspace}
\newcommand{\YELLOW}{\texttt{YELLOW}\xspace}
\newcommand{\GREEN}{\texttt{GREEN}\xspace}
\newcommand{\BLACK}{\texttt{BLACK}\xspace}

\renewcommand{\sup}[1]{\ensuremath{{}^{(#1)}}}

%% file: frontmatter.tex
% Page heads
\markboth{Aruliah {\sl et al.}}{\pampac: A Parallel Adaptive Method for Pseudo-Arclength Continuation}

% Title portion
\title{\pampac: A Parallel Adaptive Method for Pseudo-Arclength Continuation}
\author{D.A.~Aruliah
\affil{University of Ontario Institute of Technology}
Lennaert van Veen
\affil{University of Ontario Institute of Technology}
Alex Dubitski
\affil{Amadeus R\&D, Toronto}}

% Metadata Information
\acmVolume{0}
\acmNumber{0}
\acmArticle{0}
\acmYear{2012}
\acmMonth{0}

\begin{abstract}
{Pseudo-arclength continuation is a well-established method for
generating a numerical curve approximating the solution of an underdetermined system of nonlinear equations. It is an inherently sequential predictor-corrector method 
in which new approximate solutions are extrapolated from previously 
converged results and then iteratively refined.
Convergence of the iterative corrections is guaranteed only
for sufficiently small prediction steps.
In high-dimensional systems, corrector steps are extremely
costly to compute and the prediction step-length 
must be adapted carefully to avoid failed steps or 
unnecessarily slow progress.
We describe a parallel method for adapting the step-length employing several
predictor-corrector sequences of different step lengths computed 
concurrently. In addition, the algorithm permits
intermediate results of unconverged correction sequences to seed new 
predictions. This strategy results in an aggressive optimization of
the step length at the cost of redundancy in the concurrent computation.
We present two examples of convoluted solution curves of 
high-dimensional systems showing that speed-up by a factor of two
can be attained on a multi-core CPU while a factor of three is attainable on
a small cluster.}
\end{abstract}

\category{G.1.5}{Numerical Analysis}{Roots of Nonlinear Equations---Continuation (homotopy) methods}
\category{G.4}{Mathematical Software}{Parallel and Vector Implementations}

\terms{Algorithms, Performance}

\keywords{Pseudo-arclength Continuation, Parallel Computing, Adaptivity}

\acmformat{Aruliah, D. A., van Veen, L., Dubitski, A. 2012. \pampac: A Parallel Adaptive Method for Pseudo-Arclength Continuation.}

\begin{bottomstuff}
This work is supported by the Natural Sciences and Engineering Research Council.
Author's addresses: Faculty of Science, UOIT, 2000 Simcoe Street North, Oshawa, ON, L1H 7K4
\end{bottomstuff}

\maketitle

%% file: introduction.tex
\section{Introduction}

\textit{Continuation} or \textit{homotopy} problems arise naturally in numerous application domains.
They are used to study the parameter-dependence of solutions of nonlinear problems by continuously morphing between different systems of equations.
For instance, in the numerical study of the Navier-Stokes equations for fluid motion, the Reynolds number \Re
often appears as a parameter. Certain computations---\eg, those of time-periodic
solutions or travelling waves---are significantly less challenging at low Reynolds numbers than at high Reynolds numbers where dynamical
processes occur on a larger range of spatial scales.
Continuation can be used to extend the results obtained at some low Reynolds number into the physically more interesting regime by
constructing a homotopy between turbulent flows with distinct characteristics. Homotopies are also used to compute and contrast similarities and differences of flows in disparate geometries (see \citeN{kawa} and references therein).
%% rewritten

%
In mathematical terms, homotopy or continuation problems are nonlinear systems of
equations where the number of equations is one fewer than the number
of variables.
That is, given a mapping $\vec{F}:\R^{n}\times\R \rightarrow \R^n$, the equation
\begin{equation} \label{eq:base}
\vec{F}(\vec{x},\lambda)=\vec{0}
\end{equation}
with vector $\vec{x}\in\R^{n}$ and scalar $\lambda\in\R$ defines a continuation problem.
In geometric terms, Eq.~\eqref{eq:base} implicitly defines a one-dimensional
curve in $\R^{n+1}$ under suitable smoothness and consistency
properties of the mapping $\vec{F}$.
The essential idea of continuation, or homotopy, is to follow this curve of solutions in~\eqref{eq:base}.
We assume the parameter $\lambda$ lies in some specified interval $[\lambi,\lambf]\subseteq\R$.
Often, the problem can be formulated so that solving Eq.~\eqref{eq:base} for $\vec{x}$ is easy for some $\lambda^*\in[\lambi,\lambf]$ and that the goal is to obtain a solution $\vec{x}$ when $\lambda=\lambi$ or $\lambda=\lambf$.
Thus, continuation is the process of gradually morphing the solution of a
straightforward problem into the solution of a formidable problem in small parameter increments.

Numerical continuation refers to families of numerical algorithms for generating points on the solution curve.
Natural continuation and pseudo-arclength continuation are examples of predictor-corrector methods \cite{Allgower2003,Govaerts2000,Kuznetsov1998}.
In the prediction stage, a putative new point on the curve is computed and, in the correction stage, the putative candidate is iteratively
refined until a solution of~\eqref{eq:base} is found.
Such algorithms are inherently sequential: previously computed points are used to predict the next solution point on the curve.
Parallelization can be introduced within corrector iterations but the predictor steps need to be computed in sequence.

More importantly, predictor-corrector methods rely on adaptive step-size selection to make optimal progress in moving along the curve of solutions~\cite{Allgower2003}.
Strategies for adapting step-sizes selection are largely heuristic:
when corrector steps fail to converge, the predictor step is rejected,
the step-size is reduced, and a new prediction step is generated.
This turns out to be the principal bottleneck in many continuation
problems: computation time devoted to nonconvergent corrector
steps is wasted.
As such, identifying strategies to improve the performance of predictor-corrector schemes by reducing the cost of failed predictor steps is a significant challenge in modern High-Performance Computing for numerical continuation problems.

We describe in the present work a parallel software library and the underlying algorithms that extend adaptive predictor-corrector methods to amortize the cost of rejected predictor steps.
Specifically, we compute several predictor steps of different step-sizes in parallel on distinct processors. At the same time, intermediate corrector iterates can seed new predictor steps.
This strategy is most effective when corrector steps are costly.
In particular, the time for a single corrector iteration should be much larger than the communication time between processors and should not depend sensitively on the step-size.
Moreover, the curve defined by~\eqref{eq:base} should have curvature that varies dramatically so that the optimal continuation step-size also changes significantly along the solution curve.

%% file: background.tex
\section{Background}\label{sec:back}

We briefly review existing numerical continuation algorithms---notably natural parameter and pseudo-arclength continuation---before describing our parallel adaptive algorithm.
The prototypical problem is of the form \eqref{eq:base} where $\lambda\in[\lambi,\lambf]\subseteq\R$, $\vec{x}\in\R^n$ and $\vec{F}:\R^{n}\times\R \rightarrow \R^n$.
We sometimes write the underlying equation \eqref{eq:base} in the form
\begin{equation}
	\vec{F}(\vec{z})=\vec{0},\ \text{where}\ \vec{z}=(\vec{x},\lambda)\in\R^{n+1}.
\end{equation}
That is, we treat the concatenation of the $n$-vector $\vec{x}$ and the
scalar $\lambda$ as a single $(n+1)$-vector $\vec{z}$ while using the
same symbol $\vec{F}$ to denote the mapping;
technically, this notation is ambiguous but the
meaning is clear from the context.
%
%Moreover
Also, we denote appropriately-sized matrices representing Jacobian derivatives by
\begin{equation}
	\vec{F}_{\lambda} \equiv \frac{\partial\vec{F}}{\partial\lambda}\in\R^{n\times1},\quad
	\vec{F}_{\vec{x}} \equiv
        \frac{\partial\vec{F}}{\partial\vec{x}}\in\R^{n\times n},\quad\text{and}\quad
	\vec{F}_{\vec{z}} \equiv \frac{\partial\vec{F}}{\partial\vec{z}}\in\R^{n\times (n+1)}\quad\text{respectively}.
\end{equation}

For concreteness, the goal is to find a vector $\vec{x}=\vec{x}(\lambf)$ satisfying $F(\vec{x}(\lambf),\lambf)=\vec{0}$ when $\lambda$ initially starts from $\lambda^*=\lambi$ (that is, the curve is traversed with $\lambda$ increasing, at least initially).
As a first obvious method for constructing the numerical curve of solutions, one can increment the parameter $\lambda$ gradually from $\lambi$ to $\lambf$.
That is, given a point $(\vec{x},\lambda)\in\R^{n+1}$ known to lie on the curve, a new point on the curve $({\bm\xi},\lambda+h)$ is found
\begin{enumerate}
\item by incrementing $\lambda$ by a small amount $h>0$; and
\item by solving the $n$ equations $\vec{F}({\bm\xi},\lambda+h)=\vec{0}$ for the unknown ${\bm\xi}\in\R^{n}$.
\end{enumerate}
This approach is referred to as \textit{natural parameter continuation} (or, more simply, \textit{natural continuation}) \cite{Allgower2003,Govaerts2000,Kuznetsov1998}.
Algorithm~\ref{alg:npc} is a high-level description of natural continuation.
Fixing the known parameter value $\lambda\mapsto\lambda+h$ in the underdetermined nonlinear system $\vec{F}\left({\bm\xi},\lambda\right)= \vec{0}$ yields a closed $n\times n$ system in the $n$ unknown components of ${\bm\xi}$.
Intuitively, if $h$ is sufficiently small, the vector ${\bm\xi}$ should be easy to obtain using a quadratically-convergent Newton iteration (see, \eg, \citeN{Kelley2003}).
%

%%%%%%%%%%%%%%%%%%%%%%%%%%%%%%%%%%%%%%%%%%%%%%%%%%%%%%%%%%%%%%%%%%%%%%%%%%%%%%%
% Algorithm for natural parameter continuation
\IncMargin{1.5em}
\begin{algorithm}
\KwIn
{
$[\lambi,\lambf]\subseteq\R$;\ %
vector $\vec{x}\sup{0}\in\R^n$ where $\vec{F}(\vec{x}\sup{0},\lambi)=\vec{0}$;\ %
step-size $h>0$
}
\KwOut
{
vector $\vec{x}\sup{k}\in\R^n$ and scalar $\lambda\sup{k}\ge\lambf$ such that $\vec{F}\left(\vec{x}\sup{k},\lambda\sup{k}\right)=\vec{0}$
}
$k\leftarrow0$ \;%
$\lambda\sup{0}\leftarrow\lambi$
\Comment*[r]{initialization}
\While(\Comment*[f]{loop to generate successive points on curve})
{
$\lambda\sup{k}<\lambf$
}
{
$\lambda\leftarrow \lambda\sup{k}+h$ \label{ln:stepnpc}
 \Comment*[r]{predictor step} %
Iteratively solve $n\times n$ nonlinear system of equations
\[ \vec{F}\left({\bm\xi},\lambda\right)= \vec{0} \]
to obtain ${\bm\xi}\in\R^{n}$ starting from initial iterate ${\bm\xi}\sup{0}=\vec{x}\sup{k}$ \label{ln:corrnpc1} \Comment*[r]{corrector steps} 
  \eIf(\Comment*[f]{accept next point on curve})
  {iteration in line~\ref{ln:corrnpc1} converges to ${\bm\xi}$}
  {$\vec{x}\sup{k+1}\leftarrow{\bm\xi}$ \;%
   $\lambda\sup{k+1}\leftarrow\lambda$ \;%
   $k\leftarrow k+1$\;}
% Else branch
  {Reduce step-size $h$ \Comment*[f]{reject predictor step \& repeat}}
}
\Return{$\vec{x}\sup{k}$, $\lambda\sup{k}$}
\caption{Natural parameter continuation.} \label{alg:npc}
\end{algorithm}
\DecMargin{1.5em}

% Algorithm for natural parameter continuation (fin.)
%%%%%%%%%%%%%%%%%%%%%%%%%%%%%%%%%%%%%%%%%%%%%%%%%%%%%%%%%%%%%%%%%%%%%%%%%%%%%%%

Natural continuation is conceptually simple and easy to implement; however, it breaks down when the solution curve admits a \textit{fold point} (\ie, a point $(\vec{x},\lambda)$ where the Jacobian matrix $\vec{F}_{\vec{x}}(\vec{x},\lambda)$ is singular; see \citeN{Dickson2007} for an alternative characterization of fold points).
At a fold point, systematically incrementing $\lambda$ as in Algorithm~\ref{alg:npc} yields an inconsistent system of nonlinear equations that cannot be solved regardless of how small $h$ is.
Fold points do occur in practical continuation problems so other
continuation strategies need to be devised \cite{Yang1986,auto}.

\textit{Pseudo-arclength continuation} (as outlined in
Algorithm~\ref{alg:pac}; see \citet{Allgower2003,Dickson2007,Keller1977}) is a standard approach to circumvent fold point singularities.
Under the assumption that both $\vec{x}$ and $\lambda$ are smooth
functions of arclength, pseudo-arclength continuation uses
the unit direction $\vec{T}\in\R^{n+1}$ tangent to the curve for prediction.
The term ``pseudo-arclength'' applies because the step-size $h$---\ie, the Euclidean distance in $\R^{n+1}$ between successive points on the numerical curve---approximates the arclength as measured along the curve.
The tangent direction $\vec{T}$ can be determined by computing a null vector of the $n\times(n+1)$ Jacobian matrix $\vec{F}_{\vec{z}}$ or by computing finite differences between successive points on the curve \cite{Allgower2003}.
The underdetermined nonlinear system of equations \eqref{eq:base} is then closed by requiring that the solution $(\vec{x},\lambda)$ sought must lie in the hyperplane orthogonal to the tangent direction $\vec{T}$ at distance $h$ from the last known point; see line~\ref{ln:corrpac1} of Algorithm~\ref{alg:pac}.
%

%%%%%%%%%%%%%%%%%%%%%%%%%%%%%%%%%%%%%%%%%%%%%%%%%%%%%%%%%%%%%%%%%%%%%%%%%%%%%%%
% Algorithm for pseudo-arclength continuation
%
\IncMargin{1.5em}
\begin{algorithm}
\KwIn
{
$[\lambi,\lambf]\subseteq\R$;\ %
vector $\vec{x}\sup{0}\in\R^n$ where $\vec{F}(\vec{x}\sup{0},\lambi)=\vec{0}$;\ %
step-size $h>0$
}
\KwOut
{
vector $\vec{x}\sup{k}\in\R^n$ and scalar $\lambda\sup{k}\ge\lambf$ such that $\vec{F}\left(\vec{x}\sup{k},\lambda\sup{k}\right)=\vec{0}$
}
$k\leftarrow0$ \;%
$\lambda\sup{0}\leftarrow\lambi$ \Comment*[r]{initialization}
$\vec{z}\sup{0}\leftarrow\left(\vec{x}\sup{0},\lambda\sup{0}\right)\in\R^{n+1}$ \;%
Determine approximate tangent vector $\vec{T}\sup{0}\in\R^{n+1}$ at $\vec{z}\sup{0}$\;
\While(\Comment*[f]{loop to generate successive points on curve}){$\lambda\sup{k}<\lambf$}{
$\vec{w}\leftarrow \vec{z}\sup{k}+h\vec{T}\sup{k}$ \label{ln:steppac}
\Comment*[r]{predictor step}
Iteratively solve $(n+1)\times (n+1)$ nonlinear system of equations
\[ \begin{gathered} \vec{F}\left(\bm{\zeta} \right)= \vec{0}, \\
 {\vec{T}\sup{k}}^{T} \left( \bm{\zeta} - \vec{z}\sup{k} \right) = h
\end{gathered} \]
%\Indp 
to obtain $\bm{\zeta}\in\R^{n+1}$ starting from initial iterate $\bm{\zeta}\sup{0}=\vec{w}$ \label{ln:corrpac1}
\Comment*[r]{corrector steps} 
  \eIf(\Comment*[f]{accept next point on curve})
  {iteration in line~\ref{ln:corrpac1} converges to $\bm{\zeta}$}
  {
   $\vec{x}\sup{k+1}\leftarrow{\bm\zeta}_{1\!\colon\!\!n}$%
   \Comment*[r]{extract subvector from ${\bm\zeta}$} 
   $\lambda\sup{k+1}\leftarrow \zeta_{n+1}$%
   \Comment*[r]{extract last element from ${\bm\zeta}$} 
   $\vec{z}\sup{k+1}\leftarrow\left(\vec{x}\sup{k+1},\lambda\sup{k+1}\right)$ \;%
Determine approximate tangent vector $\vec{T}\sup{k+1}\in\R^{n+1}$ at $\vec{z}\sup{k+1}$\;
  Heuristically adjust the step-size $h$\;
   $k\leftarrow k+1$}
        {Reduce step-size $h$ \Comment*[f]{reject predictor step $\vec{w}$ and repeat}}
}
\Return{$\vec{z}\sup{k}=\left(\vec{x}\sup{k}, \lambda\sup{k}\right)$}
\caption{Pseudo-arclength continuation.}
\label{alg:pac}
\end{algorithm}
\DecMargin{1.5em}
%
% Algorithm for pseudo-arclength continuation (fin.)
%%%%%%%%%%%%%%%%%%%%%%%%%%%%%%%%%%%%%%%%%%%%%%%%%%%%%%%%%%%%%%%%%%%%%%%%%%%%%%%

Both natural continuation and pseudo-arclength continuation fit into a broader framework of \textit{predictor-corrector methods} \cite{Allgower2003}.
Predictor-corrector methods involve three important components:
\begin{itemize}
\item a \textit{predictor step} of a prescribed step-size;
\item a sequence of \textit{corrector steps} (alternatively \textit{corrector iterations}); and
\item an adaptive step-size selection strategy. % standard algorithms (AUTO, Matcont) always tweak step-size based on curvature, Newton convergence, etc.
\end{itemize}
The predictor step is used to seed the iterative solution of a nonlinear system of equations (successive iterates being called {corrector steps}).
There are a variety of ways in which the predictor and corrector steps can be chosen; see \citet{Allgower2003,auto}.
For instance, the predictor step in natural continuation comprises incrementing the scalar $\lambda$ by $h$; in pseudo-arclength continuation, prediction involves traversing a distance $h$ along $\vec{T}$, the direction tangent to the curve. 
Whatever the particular predictor-corrector strategy, when the
step-size $h$ is too large, the inner corrector iterations can
stagnate or diverge (either because the initial predictor is too far from a solution or perhaps because the system is inconsistent).
In either event, the predictor step is rejected and the step-size is
reduced according to a simple heuristic such as $h \leftarrow t h$
where $t$ is a user-specified parameter satisfying $0<t<1$.
Consult \citeN[Ch.~6]{Allgower2003} for more detailed strategies for adapting step-sizes.

The structure of predictor-corrector methods has to two significant consequences:
\begin{enumerate}
\item \textit{The most expensive part of the algorithm is computation
    of the corrector steps} (especially when each iteration requires the solution of a linear system of equations as in, \eg, Newton's method).
\item \textit{Rejected predictor steps are costly}
  because numerous corrector steps are computed prior to rejection.
\end{enumerate}
For moderate-sized nonlinear systems, the linear systems to be solved at each corrector iteration are amenable to dense, direct linear algebra solvers (\eg, as found in LAPACK libraries \cite{lapack});
as such, failed predictor steps may not be so punitive.
However, for large-scale systems, time (and possibly memory) requirements for each corrector iteration grows algebraically, so using direct solvers becomes infeasible.
In that case, \textit{Krylov subspace iterations} (see, \eg, \cite{saad}) can be applied within Newton iterations to determine the corrector steps within a continuation algorithm.
This process is referred to as \textit{Newton-Krylov continuation} \cite{sanchez,keyes}). The convergence of Krylov subspace methods depends
critically on the properties of the Jacobian matrix and may require preconditioning.
For certain problems with upward from $10,000$ degrees of freedom, individual corrector steps in systems can require hours (in some cases, days) of computation---even using a well-tuned Krylov subspace method---and the penalty incurred for failed predictor steps is prohibitive.
Given that rejected predictor steps can result from using a large step-size, an obvious strategy is to use very small step-sizes.
Unfortunately, using small step-sizes impedes progress along the curve from $\lambi$ to $\lambf$, both in natural and in pseudo-arclength continuation.
Thus, efficient step-size adaptivity requires trading off between these conflicting concerns.

To achieve this balance, we develop a parallel software library---\pampac, a Parallel Adaptive Method for Pseudo-Arclength Continuation{}---that permits adaptive step-size selection within the predictor-corrector framework of pseudo-arclength continuation.
The parallel algorithm implemented can be used for different kinds of correction steps and leads to speed-up for problems in which the computation of a correction iteration takes much longer than communication between processors.
It is particularly effective for continuation problems in which the solution curve exhibits large variations in curvature (\ie, where the optimal step-size changes wildly along the curve).

%% file: parallelization-scheme.tex
\section{Scheme for parallelization}\label{sec:parallelization}

There are two essential ways to achieve parallelism in the selection of step-sizes for predictor steps.
The first is most obvious: use a concurrent sequence of predictor steps of
different step-sizes $t_{1}h$, $t_{2}h$, \ldots, $t_{W}h$ for some
prescribed positive scalars $\{t_{\alpha}\}_{\alpha=1}^{W}$.
That is, given $W$ processes, an initial point $\vec{z}\in\R^{n+1}$ known
to satisfy $\vec{F}(\vec{z})=\vec{0}$ up to some fixed tolerance, and a unit tangent direction
$\vec{T}\in\R^{n+1}$, each process $\alpha$ computes a predictor step
${\bm\zeta}_{\alpha}^{(0)}=\vec{z} + t_{\alpha} h\vec{T}$.
A maximum step-size $h_{\max}$ may have to be set to avoid
spurious convergence to a remote branch of the continuation curve;
in that instance, fewer than $W$ predictor steps would be computed.
Once a prediction ${\bm\zeta}_{\alpha}^{(0)}$ is determined by process $\alpha$,
the process computes a sequence of corrector steps
$\{ {\bm\zeta}_{\alpha}^{(0)}, {\bm\zeta}_{\alpha}^{(1)},
{\bm\zeta}_{\alpha}^{(2)},\dots\}$.
Subsequent corrector iterations do not require inter-process communications and thus can be carried out concurrently by distinct processes.
Each process $\alpha$ maintains its own iteration counter $\nu_\alpha$ as
well as the iterates ${\bm\zeta}_{\alpha}^{(\nu_\alpha)}$ and the associated nonlinear residuals $\vec{r}_{\alpha}^{(\nu_\alpha)}\equiv
\vec{F}\left({\bm\zeta}_{\alpha}^{(\nu_\alpha)}\right)$.
A second strategy for parallel step-size adaptivity is to
use intermediate computations to seed new predictions.
That is, suppose process $\alpha$ is initialized using the predictor
${\bm\zeta}_{\alpha}^{(0)}=\vec{z}_{\alpha} + h_{\alpha} \vec{T}_{\alpha}$
and proceeds to compute a sequence of correction steps $\{
{\bm\zeta}_{\alpha}^{(0)}, {\bm\zeta}_{\alpha}^{(1)},
{\bm\zeta}_{\alpha}^{(2)},\dots\}$.
Rather than waiting for the corrector iterations to converge, the
intermediate iterates can be used to compute a normalized secant direction
$\hat{\vec{T}}_{\alpha}$ in the direction of
${\bm\zeta}_{\alpha}^{(\nu_\alpha)}-\vec{z}_{\alpha}$
for some iterate $\nu_\alpha$.
Assuming the corrector iterates are sufficiently close to the curve of
solutions, the new secant direction $\hat{\vec{T}}_{\alpha}$
approximates a tangent to the curve and can be handed off to another
process $\beta$ to generate a new 
predictor point $\vec{z}_{\beta} +  h_{\beta}
{\vec{T}}_{\beta}$, where $\vec{z}_{\beta}={\bm\zeta}_{\alpha}^{(\nu_\alpha)}$,
$h_\beta=t_{\ell}h_\alpha$ for some $\ell\in\{1,\dotsc,W\}$, and
$\vec{T}_{\beta} = \hat{\vec{T}}_{\alpha}$.
To manage concurrent processes when both these strategies are employed, we represent each process by a node in a rooted tree with a master process at the root.
The user specifies the tree's width $W>0$ and depth $D>0$;
the width corresponds to the number of scalars $t_{1}$, $t_{2}$, \ldots, $t_{W}$ that multiply the step-length $h$ and the depth corresponds to the number of extrapolated predictor steps computed from intermediate corrector iterates.
New nodes, {\ie}, new corrector sequences, are seeded only if
physical processors are available; the queueing of multiple processes
on a single processor hinders the continuation.
Each node in the tree corresponds to a potential computational process.
A node $\alpha$ is associated with a current iteration counter $\nu_\alpha$ and a current corrector iterate ${\bm\zeta}_{\alpha}^{(\nu_\alpha)}\in\R^{n+1}$.
The node $\alpha$ is generated with
\begin{equation}
\nu_\alpha\leftarrow 0 \text{\ and\ } {\bm\zeta}_{\alpha}^{(0)} \leftarrow
\vec{z}_{\alpha}^{\text{init}} + h_{\alpha}^{\text{init}} \vec{T}_{\alpha }^{\text{init}},
\end{equation}
where the data $\vec{z}_{\alpha}^{\text{init}}$, $\vec{T}_{\alpha}^{\text{init}}$, and ${h}_{\alpha}^{\text{init}}$ are determined from the parent node (\ie, the previous point on or near the continuation curve).
When created, the node $\alpha$ also records
$\nu_{\alpha}^{\text{init}}$---the number of corrector iterations that its parent had computed prior to spawning node $\alpha$.

Node $\alpha$ also maintains a real parameter $h_\alpha$ representing the base step-size node that $\alpha$ uses to construct new predictor steps.
Initially, $h_\alpha\leftarrow h_{\alpha}^{\text{init}}$, \ie, the base step-size node $\alpha$ uses to make predictor steps is the same as the step-size used to initialize node $\alpha$.
However, in the event that all child nodes spawned from $\alpha$ lead to divergent sequences, $h_\alpha$ is reduced by some constant scaling factor $t\in(0,1)$ before spawning more predictions.
The scaling factor $t$ is fixed to ensure all new prediction steps are at a distance shorter than the shortest of the steps that just failed:
\[ t = 0.9\, \frac{t_{\min}}{t_{\max}},
\ \text{where}\ t_{\min}=\min_{1\le k\le W}  t_{k}
\ \text{and}\ t_{\max}=\max_{1\le k\le W}  t_{k}. \]
This precaution ensures new predictors generated do not repeat earlier computations (after the failed child nodes are all deleted).
If the base step-size $h_{\alpha_{r}}$ of the root node is reduced below some user-specified threshold $\mathtt{H\_MIN}$, the continuation algorithm halts.

Each node $\alpha$ is assigned a color $c_\alpha \in \left\{ {\GREEN},
{\YELLOW},
{\RED},
{\BLACK} \right\}$ as the algorithm proceeds.
The color $c_\alpha$ of node $\alpha$ is determined using
 {$\norm{\vec{r}_{\alpha}^{(\nu_\alpha)}}_{2}$}, 
\ie,  the 2-norm of the nonlinear residual
$\vec{r}_{\alpha}^{(\nu_\alpha)} = \vec{F}\left({\bm\zeta}_{\alpha}^{(\nu_\alpha)}\right)\in \R^{n}$ at the current corrector iterate ${\bm\zeta}_{\alpha}^{(\nu_\alpha)}$.
Given positive parameters $\mathtt{TOL\_RESIDUAL}$, ${\mathtt{GAMMA}}$, and ${\mathtt{MU}}$, the rules are as follows:
\begin{enumerate}
\item $c_\alpha\leftarrow{\GREEN}$ if {$\norm{\vec{r}_{\alpha}^{(\nu_\alpha)}}_{2}\le\mathtt{TOL\_RESIDUAL}$} (\ie, the current iterate
  ${\bm\zeta}_{\alpha}^{(\nu_\alpha)}$ is deemed to have converged to a point on or sufficiently near the solution curve);
\item \label{i:yellow} $c_\alpha\leftarrow{\YELLOW}$ if {$\norm{\vec{r}_{\alpha}^{(\nu_\alpha)}}_{2}^{\mathtt{GAMMA}} \le\mathtt{TOL\_RESIDUAL}$} (\ie, the \textit{next iterate} ${\bm\zeta}_{\alpha}^{(\nu_\alpha+1)}$ is expected to have converged);
\item \label{i:black} $c_{\alpha}\leftarrow{\BLACK}$ if $\nu_{\alpha}>\mathtt{MAX\_ITER}$ or {$\norm{\vec{r}_{\alpha}^{(\nu_\alpha)}}_{2} > {\mathtt{MU}} \norm{\vec{r}_{\alpha}^{(\nu_\alpha-1)}}_{2}$}
(\ie, the maximum number of corrector iterations is exceeded or the reduction of the residual is insufficient in consecutive corrector iterations); and
\item $c_{\alpha}\leftarrow{\RED}$ otherwise.
\end{enumerate}
The efficiency of the parallelization depends on the preceding criteria for coloring nodes.

New nodes are colored \RED by default;
the colors are reassessed after the current corrector iterates and the corresponding nonlinear residuals are computed on every node.
The specific criteria for coloring a node \GREEN, \YELLOW, or \BLACK depend on the nature
of the corrector steps and the continuation problem at hand.
As elucidated in Sec.~\ref{sec:structure}, the user chooses these criteria
by providing problem-dependent parameters $\mathtt{TOL\_RESIDUAL}$, $\mathtt{MAX\_ITER}$, ${\mathtt{GAMMA}}$, and ${\mathtt{MU}}$.
The parameters $\mathtt{TOL\_RESIDUAL}$ and $\mathtt{MAX\_ITER}$ are standard as expected in any iterative solver.
The other parameters ${\mathtt{GAMMA}}$ and ${\mathtt{MU}}$ are based on the asymptotic behavior of the user's choice of corrector sequences.
For instance, when using Newton's method for corrector steps,
we expect the decrease in the residual to be quadratic (or at least superlinear)
near the solution curve; in that case, it makes sense to choose ${\mathtt{GAMMA}}\in(1,2)$ in the criterion~\eqref{i:yellow} to color ``nearly converged nodes'' {\YELLOW}.
For predictor steps too far from the solution curve,
we often see linear convergence only;
it makes sense, then, to choose $\mathtt{MU}\in(0,1)$ in the criterion~\eqref{i:black} to color a node {\BLACK} when it is not converging sufficiently quickly.
Nodes colored \BLACK get deleted from the tree since computing more corrector steps likely leads to divergence or spurious convergence on another branch of the continuation curve.
%
%This is where problem-dependent knowledge helps optimize parallel speed-up.
%
%The numerical examples we present in Section~\ref{sec:num} use Newton and Newton-Krylov iteration, with $\mathtt{GAMMA}=2$ and $\mathtt{MU}=1/2$.
%

%
Algorithm~\ref{alg:pampac} is a high-level summary of the essential
steps underlying our parallel adaptive algorithm.
The main conceptual pieces are the concurrent computation of corrector steps from
predictor steps of various step-lengths, extrapolation from intermediate corrector steps,
and management of the parallel computations using a rooted tree with colored nodes.
%%%%%%%%%%%%%%%%%%%%%%%%%%%%%%%%%%%%%%%%%%%%%%%%%%%%%%%%%%%%%%%%%%%%%%%%%%%%%%%
% Algorithm for parallel adaptive continuation
\IncMargin{1.5em}
\begin{algorithm}
\KwIn
{
$[\lambi,\lambf]\subseteq\R$;\ %
vector $\vec{x}^{*}\in\R^{n}$ where
$\vec{F}(\vec{x}^{*},\lambi)=\vec{0}$;\ %
step-size $h\neq0$;\ %
tangent direction $\vec{T}^{*}\in\R^{n+1}$ where
$\norm{\vec{T}^{*}}_2=1$;\ %
tolerance $\mathtt{TOL\_RESIDUAL}>0$;\ %
$\mathtt{GAMMA}>0$;\ %
$\mathtt{MU}>0$;\ %
depth $D>0$;\ %
width $W>0$;\ %
positive scaling parameters
$\{t_{1},t_{2},\dotsc,t_{W}\}$
}
\KwOut
{
vector $\vec{x}\sup{k}\in\R^n$ and scalar $\lambda\sup{k}\ge\lambf$ such that $\vec{F}\left(\vec{x}\sup{k},\lambda\sup{k}\right)=\vec{0}$
}
\label{ln:seed}
Seed root node $\alpha_{r}$ with data 
$(\vec{z}_{\alpha_r},\hat{\vec{T}}_{\alpha_r},h_{\alpha_r}) \leftarrow
((\vec{x}^{*},\lambi),\vec{T}^{*},h)$
\Comment*[r]{initialization}%
\Repeat(\Comment*[f]{loop to generate successive points on curve})
{root node $\alpha_r$ has $\lambda\sup{\alpha_r}\ge\lambf$}
{
\ForEach(\Comment*[f]{Spawn new nodes}){leaf node $\alpha$}
{%
\label{ln:forstart}%
\If{depth of $\alpha < D$}
{
% added: \zeta does not exist on alpha_r
\If{$\alpha\neq \alpha_r$}{
$\hat{\vec{T}}_{\alpha}^{\text{init}} \leftarrow
\norm{{\bm\zeta}_{\alpha}^{(\nu_\alpha)}-\vec{z}_{\alpha}^{\text{init}}}_2^{-1}
\left( {\bm\zeta}_{\alpha}^{(\nu_\alpha)}-\vec{z}_{\alpha}^{\text{init}} \right)$
\Comment*[r]{secant direction}}
\For{$\ell=1\,\colon W$}
{
\If{processors are available}
{
Assign processor $\beta_{\ell}$ the data
\[\begin{gathered}
  \vec{z}_{\beta_\ell}^{\text{init}} \leftarrow {\bm\zeta}_{\alpha}^{\nu_\alpha},\quad
  \hat{\vec{T}}_{\beta_\ell}^{\text{init}} \leftarrow \hat{\vec{T}}_{\alpha},\quad
  h_{\beta_\ell}^{\text{init}} \leftarrow t_{\ell} h_{\alpha},\quad
  \nu_{\beta_\ell}^{\text{init}} \leftarrow \nu_\alpha, \\
    \nu_{\beta_\ell} \leftarrow 0,\quad
  {\bm\zeta}_{\beta_\ell}^{(0)} \leftarrow 
  {\bm\zeta}_{\alpha}^{\nu_\alpha} + t_{\ell} h_{\alpha} \hat{\vec{T}}_{\alpha},\quad
  c_{\beta_\ell} \leftarrow \RED
\end{gathered}\]
\label{ln:forend}} % end of "If processors" block
} % end of "For \ell" block
} % end of "If depth" block
} % end of "ForEach leaf" block
\label{ln:corr}%
Compute single corrector steps concurrently on all \RED
and {\YELLOW} nodes of the tree. \;
\label{ln:update}%
Traverse nodes of tree updating residuals
\[ \vec{r}_{\alpha}^{(\nu_{\alpha})} \leftarrow
  \vec{F}\left({\bm\zeta}_{\alpha}^{(\nu_{\alpha})}\right)\]
and updating colors $c_\alpha$ accordingly. \;
$\text{PruneTree}(\alpha_{r})$, deleting \BLACK nodes and eliminating redundant subtrees.\;
\label{ln:prune}%
\While(\Comment*[f]{update root})
{root node $\alpha_r$ has a single child, $\beta$, and $\beta$ is \GREEN}
{
Write the solution on $\alpha_r$ to disk\;
Update root: $\alpha_r \leftarrow \beta$
}
} % end of "Repeat Until" block
\caption{Essential parallel adaptive algorithm} \label{alg:pampac}
\end{algorithm}
\DecMargin{1.5em}
% Algorithm for parallel adaptive continuation (fin.)
%%%%%%%%%%%%%%%%%%%%%%%%%%%%%%%%%%%%%%%%%%%%%%%%%%%%%%%%%%%%%%%%%%%%%%%%%%%%%%%

%
The outermost loop of Algorithm~\ref{alg:pampac} corresponds to accepting and logging valid points along the continuation curve.
The nested loop spanning lines~\ref{ln:forstart} to~\ref{ln:forend} describe the generation of new predictor points and their distribution over available processors.
The processors compute individual corrector steps concurrently in line~\ref{ln:corr} requiring synchronization before advancing to line~\ref{ln:update}.
This implies that the parallelization scheme is efficient only when the time required for a corrector step is much greater than the communication costs and when the time for the concurrent corrector computations is roughly the same on all processors.
At a high level, the algorithm is mostly straightforward;
further explanation is required to understand the pruning algorithm in line~\ref{ln:prune}.

Concurrent computations are managed by pruning the tree in two stages.
In the first stage, all subtrees rooted at \BLACK nodes are removed.
The corrector sequences associated with \BLACK nodes are deemed to have made insufficient progress;
as such, computing more corrector steps likely leads to divergence or spurious convergence to a point on another branch of the continuation curve.
Other non-\BLACK nodes are unaltered in this first stage.
% unless such a node has only \BLACK child nodes (\ie, only nonconvergent child corrector sequences).
%
%

%
\begin{figure}
\includegraphics[width=.9\textwidth]{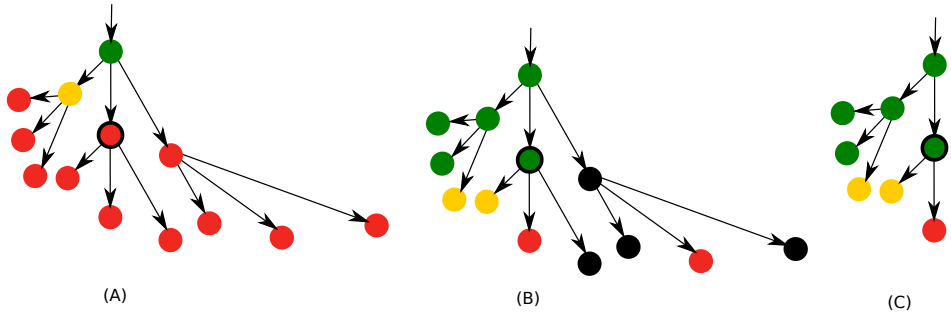}
\caption{Illustration of the first stage of the pruning algorithm on a tree of width 3 and depth 3. The links are drawn from left to
right in the order of increasing step length. (A) \GREEN root node representing a valid solution point and two levels of child nodes.
Each of the nodes at depth 2 has three \RED child nodes that have not done computed and correctoer iterations as of yet.
In (B), concurrent corrector steps have been computed on all \YELLOW and
\RED nodes. After computing and assessing the residuals on all of the nodes, some have been colored \GREEN (converged), \YELLOW (almost converged), \BLACK (diverged), or \RED (undecided).
In (C), all the \BLACK nodes and their attached subtrees have been deleted.}
\label{fig:prune_1}
\end{figure}
Figure~\ref{fig:prune_1} provides a graphical illustration of the first stage of pruning.
In panel (A), the \GREEN root node represents a valid solution point with three child nodes (one \YELLOW and two \RED nodes).
In panel (B), a concurrent corrector step is computed for all \YELLOW and \RED nodes (requiring 12 active processes in this case).
At the end of the corrector step, four nodes are now \GREEN, two are \YELLOW, and four are \BLACK.
The first stage of pruning is simply to
delete all \BLACK nodes and their associated subtrees which leads to the configuration shown in panel (C).

The second stage of pruning---line~\ref{ln:prune} of Algorithm~\ref{alg:pampac}---involves comparing child nodes of a given parent node with viable corrector iterates and deciding which to keep and which to delete (along with associated subtrees).
Doing so permits allocation of computing resources to nodes that are deemed to yield the greatest benefit.
The greatest gain is determined by balancing the longest distance travelled along the continuation curve (\ie, pseudoarclength) against 
the least amount of computational work (measured in corrector iterations).
\RED child nodes are kept by default---it is unclear whether they will yield convergent corrector sequences or not.
When deciding which \GREEN or \YELLOW child nodes to keep, a more sophisticated criterion is applied that requires a few definitions.

\begin{definition}
A {\em path} $P$ in the computation tree is a connected subtree in which each node has at most one child node.
\end{definition}
A path in the computation tree corresponds to a putative segment of the continuation curve.
An ordered sequence of connected nodes (\eg, $P=(\alpha_1,\alpha_2,\dotsc,\alpha_N)$) represents a path in the computation tree.

Being able to identify paths in the computation tree, there are two important measures needed to control our algorithm.
\begin{definition}
The {\em initialization length} of the path $P=(\alpha_1,\alpha_2,\dotsc,\alpha_N)$ is
\begin{equation}
L(P) = \sum_{k=1}^{N} h_{\alpha_k}^{\text{init}}.
\end{equation}
\end{definition}
\begin{definition}
The {\em iteration cost} of the path $P=(\alpha_1,\alpha_2,\dotsc,\alpha_N)$ is
\begin{equation}
I(P) = \mu_{\alpha_{1}}\text{\ where\ } \mu_{\alpha_{k}}=
\begin{cases}
\nu_{\alpha_k}, & \text{if\ }k=N, \\
\max\left(\nu_{\alpha_k}, \mu_{\alpha_{k+1}} + \nu_{\alpha_{k+1}}^{\text{init}}\right), & \text{otherwise}.
\end{cases}
\end{equation}
\end{definition}

That is, $L(P)$ is the pseudo-arclength of the segment of the curve associated with the nodes of $P$ in sequence---assuming the points associated with the nodes all lie on the continuation curve.
Similarly, $I(P)$ is the accumulated number of corrector
steps computed to attain the current state of the entire path.

It is useful to identify two specific kinds of paths in the tree
to choose between ``converged'' and ``nearly converged'' corrector sequences (represented on the tree by \GREEN and \YELLOW nodes, respectively).
\begin{definition}
A {\em {valid} path} is a path in which all nodes are \GREEN.
\end{definition}
\begin{definition}
A {\em {viable} path} is a path in which all nodes are either \GREEN or \YELLOW.
\end{definition}
Notice a {valid} path is, by definition, a {viable} path also.
Thus, before the second stage of pruning, {valid} and {viable} paths are identified recursively (line~\ref{ln:prune} of Algorithm~\ref{alg:pampac}).
Naturally, a {valid} path comprises nodes associated with points known to be on the continuation curve (within the desired tolerance).
{Viable} paths consist of nodes associated with points that are either known to be on the continuation curve or points that are likely to be on the curve in the next concurrent corrector iteration.
As such, when comparing a node's \GREEN and \YELLOW child nodes,
Algorithm~\ref{alg:paths} chooses the node associated with the longest {viable} or {valid} path (measured in pseudo-arclength).

In the event that distinct {valid} and {viable} paths exist,
Algorithm~\ref{alg:paths} decides whether to keep the best {valid} or the best {viable} path.
In particular, to choose between the longest valid path $P_{\text{valid}}(\alpha)$ and the longest viable path $P_{\text{viable}}(\alpha)$ extending below node $\alpha$, the rate of progress along the curve is estimated by computing
\begin{equation}\label{eq:pathchoose}
\frac{L(P_{\text{valid}}(\alpha))}{I(P_{\text{valid}}(\alpha))}%
\quad\text{and}\quad%
\frac{L(P_{\text{viable}}(\alpha))}{I(P_{\text{viable}}(\alpha))+1}.%
\end{equation}
The algorithm keeps whichever child node leads to the path with the greatest rate of progress as computed in~\eqref{eq:pathchoose}; the other \GREEN and \YELLOW child nodes are deleted (with their associated subtrees).

% Example with illustrations
\begin{figure}
\begin{center}
\includegraphics[width=.6\textwidth]{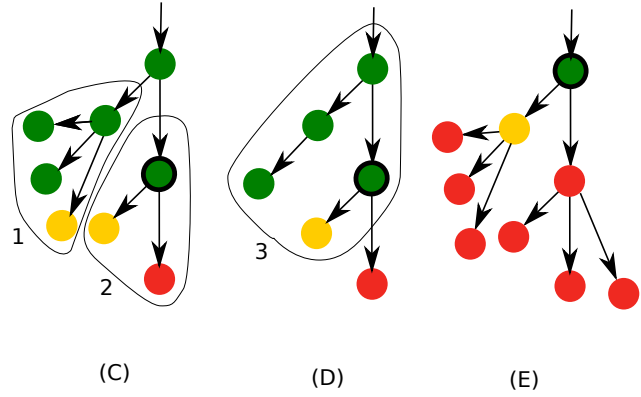}
\end{center}
\caption{Illustration of the pruning algorithm on a tree of width 3 and depth 3. For each subtree, the path that is expected to
yield the fastest progress is retained with any \RED nodes.
This path is first selected in the smaller subtrees 1 and 2 and then in the larger subtree 3.
After pruning, the \GREEN node at the root of subtree 2 in panel (C) becomes the new root node and three new leafs are spawned from both of its children (as shown in panel (E)).}
\label{fig:prune_2}
\end{figure}

Figure~\ref{fig:prune_2} illustrates Algorithm~\ref{alg:paths} for a tree of width 3 and depth 3 continuing where Figure~\ref{fig:prune_1} stopped.
The step lengths in this example are $\{t_k\}_{k=1}^W=\{1/4,1,3/2\}$.
In panel (C), the root \GREEN node represents a solution that has been computed after two corrector steps, while
its child nodes have been spawned after one step.
First, subtrees 1 and 2 are examined to choose which child to keep.
Subtree 2 in panel~(C) has only a single viable path, so it is left alone.
Subtree 1 in panel~(C) has two valid and three viable paths from its root \GREEN node.
The longest valid and longest viable paths  are compared using~\eqref{eq:pathchoose} and only the path with the greatest ``speed'' is kept.
The longest valid path took two corrector steps to construct and gives a gain of $t_1 h+t_2 t_1 h=h/2$ in arc length,
while the longest viable path is expected to yield a gain of $t_1 h+t_3 t_1 h=5h/8$ in three steps.
Since $(h/2)/2>(5h/8)/3$, the longest valid path is retained.
The result in panel~(D) is subtree~3 hanging off the root \GREEN node which has one valid and two viable paths.
 Here, retaining the all-green (valid) subtree, on the left, would lead to a gain of 
$(1+t_1+t_2 t_1)h=3h/2$ in three corrector steps. Retaining the \GREEN and \YELLOW (viable but not valid) subtree, on the right, should give
a gain of $(1+t_2+t_1 t_2)h=9h/4$ in four corrector steps. Since $(3h/2)/3=h/2< 9h/16=(9h/4)/4$, the 
the \GREEN and \YELLOW sub tree is retained.
Subsequently, the next \GREEN node on the retained path becomes the new root node in panel~(E).
The \YELLOW and \RED nodes from  panel~(D) are used to compute new predictor points in panel~(E).

%%%%%%%%%%%%%%%%%%%%%%%%%%%%%%%%%%%%%%%%%%%%%%%%%%%%%%%%%%%%%%%%%%%%%%%%%%%%%%%
% Second part of the pruning: finding the optimal path.
\IncMargin{1.5em}
\begin{algorithm}
\KwIn
{
root node $\alpha_r$ of computation tree
}
\ForEach{node $\alpha$ in a depth-first traversal of computation tree rooted at $\alpha_{r}$}
{
\If{$c_{\alpha}=\mbox{\BLACK}$}{delete subtree rooted at $\alpha$}
}
\ForEach{node $\alpha$ in a depth-first traversal of computation tree rooted at $\alpha_{r}$}
{
compute paths $P_{\text{viable}}(\alpha)$ \& $P_{\text{valid}}(\alpha)$,
the respective viable and valid paths of longest initialization path length rooted at node $\alpha$
}
\ForEach{node $\alpha$ in a depth-first traversal of computation tree rooted at $\alpha_{r}$}
{
$P_{\text{best}}(\alpha) \leftarrow \emptyset$ \newline
\nl \If{$P_{\text{viable}}(\alpha)\neq\emptyset$}
{$P_{\text{best}}(\alpha) \leftarrow P_{\text{viable}}(\alpha)$ \newline
\nl \If{$P_{\text{valid}}(\alpha)\neq\emptyset$ \and
$P_{\text{valid}}(\alpha)\neq P_{\text{viable}}(\alpha)$}
{
$S_{\text{valid}}(\alpha)\leftarrow L(P_{\text{valid}}(\alpha)) / I(P_{\text{valid}}(\alpha))$ \newline \nl
$S_{\text{viable}}(\alpha)\leftarrow L(P_{\text{viable}}(\alpha)) / \left[I(P_{\text{viable}}(\alpha))+1\right]$ \newline \nl
\If{$S_{\text{valid}}(\alpha)\ge S_{\text{viable}}(\alpha)$}
{
$P_{\text{best}}(\alpha) \leftarrow P_{\text{valid}}(\alpha)$
}
}
}
\ForEach{child node $\beta$ of node $\alpha$}
{
\If{$\beta\not\in P_{\text{best}}(\beta)$ and $c_{\beta}\in\{\YELLOW,\GREEN\}$}
{
delete subtree rooted at $\beta$
}
}
}
\Return{}
\caption{\label{alg:paths}PruneTree($\alpha_r$)}
\end{algorithm}
\DecMargin{1.5em}
% Second part of the pruning: finding the optimal path. (fin.)
%%%%%%%%%%%%%%%%%%%%%%%%%%%%%%%%%%%%%%%%%%%%%%%%%%%%%%%%%%%%%%%%%%%%%%%%%%%%%%%

%% file: structure-of-the-code.tex
\section{Structure of the code}
\label{sec:structure}

We have implemented the algorithms from Section~\ref{sec:parallelization} in a library called \pampac (Parallel Adaptive Method for Pseudo-Arclength Continuation).
The \pampac library is written in C (ISO/IEC 9899:1999) using MPI-2 libraries (the Message Passing Interface, see \cite{MPI}) for parallelization.
The main \pampac library comes with template \texttt{Makefiles} for easy building; some configuration of the file \texttt{Makefile.in} may be required to ensure that all library dependencies are met on a user's system, but the goal is straightforward and builds on any POSIX system.
The code is packaged with an example (described more fully in Section~\ref{ex:modified_ks}) to help users modify their own continuation codes for use with \pampac easily.

To use the \pampac library, the user needs to supply:
\begin{itemize}
\item a \texttt{main} function (usually in a file \texttt{main.c}) to drive the computation;
\item a C function \verb^compute_residual^ that evaluates $\vec{F}(\vec{z})$ in~\eqref{eq:base};
\item a C function \verb^single_corrector_step^ that computes the updated corrector iterate ${\bm\zeta}_{\alpha}^{(\nu_{\alpha}+1)}$ from ${\bm\zeta}_{\alpha}^{(\nu_{\alpha})}$ and $\vec{T}$ as required in line~\ref{ln:corrpac1} of Algorithm~\ref{alg:pac} for pseudo-arclength continuation; and
\item an input file (ASCII text) containing the initial data point $\vec{z}^{*}=(\vec{x}^{*},\lambda^{*})$ on the continuation curve.
\item a parameter file (usually \texttt{parameters.txt}) for tuning the behavior of Algorithm~\ref{alg:pampac}.
\end{itemize}
The subdirectory of \texttt{examples} contains a template file \texttt{main.c} that can be used with little or no modification.
That directory also contains a sample \texttt{parameters.txt} file for controlling the algorithm behavior.
The code is modular to aid in debugging and understanding; each of the core tasks listed in Algorithm~\ref{alg:pampac} of Section~\ref{sec:parallelization} is performed by a separate function.

The \texttt{main} function performs two primary tasks: it initializes and finalizes MPI communication and it divides work between the master and the slave processes (this is typical in the Single-Program-Multiple-Process (SPMP) paradigm).
The master processor parses the user-provided parameter file to determine the algorithm-tuning parameters and calls the routine \verb^master_process^;
this routine does some preprocessing (\eg, loading the initial point from the user's input file and computing an initial tangent direction) before initiating Algorithm~\ref{alg:pampac}.
The slave processes all call the function \verb^slave_process^ in which they idle until receiving data from the master process---the data being a point ${\bm\zeta}\in\R^{n+1}$ and some tangent direction $\vec{T}\in\R^{n+1}$ from which a corrector iteration can be computed.
The only interprocess communication during the continuation loop consists
of the root process sending these data to the slaves and each slave returning the
result of a corrector step to the root process.
The routines \verb^master_process^ and \verb^slave_process^ package the core components of Algorithm~\ref{alg:pampac} in a manner that alleviates the burden of managing the parallel computation from the user.
As mentioned in Section~\ref{sec:parallelization}, the corrector iterations are independent and are generally much more expensive than the cost of inter-process communication.

The user-supplied routines \verb^compute_residual^ and \verb^single_corrector_step^ have the following C function prototypes:
\begin{verbatim}
  void compute_residual (int N_dim, const double* z, double* residual);
  void single_corrector_step (int N_dim, double* zeta, double* Tangent);
\end{verbatim}
Both functions do not have return values; rather, the ``output'' values computed are passed by reference (a common idiom in C and FORTRAN programming).
Notice that, relative to the mathematical description of the template continuation problem in Section~\ref{sec:back}, \verb^N_dim^$\,= n+1$, i.e., \verb^N_dim^ refers to the dimension of the vector $\vec{z}=(\vec{x},\lambda)$ rather than the dimension of the vector $\vec{x}$.
Thus, in \verb^compute_residual^, the ``input'' values are the (integer) dimension \verb^N_dim^ of the problem and the \verb^N_dim^-vector pointed to by the pointer \verb^z^;
the ``output'' is the residual vector $\vec{F}(\vec{z})$ in~\eqref{eq:base} stored in an
array of length \verb^N_dim-1^ in memory pointed to by 
the pointer \verb^residual^.
Similarly, in \verb^single_corrector_step^, the ``input'' values are the (integer) dimension \verb^N_dim^ of the problem, the \verb^N_dim^-vector pointed to by the pointer \verb^zeta^, and the \verb^N_dim^-vector pointed to by the pointer \verb^Tangent^ (corresponding to $\vec{T}$).
After calling  \verb^single_corrector_step^, the array pointed to by \verb^zeta^ has been overwritten with the updated corrector iterate as in Algorithm~\ref{alg:pac}.
These functions need to be compiled with \texttt{main.c}---and any user-required dependencies---to produce an executable that can be run in parallel on numerous processors.
Assuming that the user's \texttt{Makefile} is suitably configured, the user can link the executable with external library functions required by their routines.

Parameters controlling the parallel continuation Algorithm~\ref{alg:pampac} are loaded
from a plain text file at run-time by the master processor.
The user needs to provide the following values:
\begin{itemize}
\item \verb+N_DIM+: the number of unknowns $n+1$;
\item \verb+LAMBDA_MIN+ and \verb+LAMBDA_MAX+: bounds on the interval $[\lambi,\lambf]$ in which the continuation parameter $\lambda$ lies;
\item \verb+LAMBDA_INDEX+: integer between \texttt{0} and \verb+N_DIM-1+ that is the index of the parameter $\lambda$ in any \verb+N_DIM+-vector;
\item \verb+DELTA_LAMBDA+: parameter for initial corrector iterations to generate a second point on the curve from the first (required to bootstrap the algorithm);
\item \verb+H_MIN+, \verb+H_MAX+ and \verb+H_INIT+: the minimal, maximal and initial pseudo-arclength step-size;
\item \verb+MAX_ITER+: the maximum number of corrector steps before coloring a node \BLACK;
\item \verb+TOL_RESIDUAL+: the threshold residual tolerance in Eq.~\ref{eq:base} for accepting \GREEN nodes (\ie, when  \smash{$\norm{\vec{r}_{\alpha}^{(\nu_\alpha)}}_{2} \le\,$}\verb+TOL_RESIDUAL+);
\item \verb+MU+: the threshold reduction in residual for \BLACK nodes (\ie, when \smash{$\norm{\vec{r}_{\alpha}^{(\nu_\alpha)}}_{2} > \mathtt{MU} \norm{\vec{r}_{\alpha}^{(\nu_\alpha-1)}}_{2}$});
\item \verb+GAMMA+: the threshold rate of residual reduction for \YELLOW nodes (\ie,
when \smash{$\mathtt{GAMMA}\log\norm{\vec{r}_{\alpha}^{(\nu_\alpha)}}_2\le\log\,$}\verb+TOL_RESIDUAL+;
\item \verb+MAX_DEPTH+: the maximum depth of the tree, $D$;
\item \verb+MAX_CHILDREN+: the maximum width of the tree, $W$;
\item \verb+SCALE_PROCESS_K+ ($\mathtt{K}=0\ldots W-1$): the step-size scalings $t_{\mathtt{K}}$ in Algorithm~\ref{alg:pampac}; and
\item \verb+VERBOSE+: an integer parameter controlling verbose output.
\end{itemize}

Certain parameters in the user's parameter file are not mentioned in the description of Algorithm~\ref{alg:pampac} from Section~\ref{sec:parallelization}.
To circumvent a stagnating loop, the user can specify a positive integer \verb^MAX_ITER^ to terminate corrector iterations.
The parameter \verb+LAMBDA_INDEX+ provides additional flexibility by permitting the user to specify any integer index of $\vec{z}$---using 0-based indexing as is conventional in C---for the continuation parameter.
That is, the parameter $\lambda$ does not need to be the $(n+1)^{\mathrm{st}}$ component of the $(n+1)$-vector $\vec{z}$.

To bootstrap Algorithm~\ref{alg:pampac},
the master processor requires an initial tangent direction in addition to  the initial point loaded from the user's input file.
It generates an approximate tangent direction by carrying out a few corrector iterations to generate another point near the initial point and computing a secant direction between those two points.
At any given point on the continuation curve, there are two anti-parallel tangent directions;
as such, the sign of  \verb+H_INIT+ is used to fix the initial direction of the continuation
(\ie, the tangent direction used to generate the second point on the curve is oriented in the direction of $\lambda$ increasing or decreasing when \verb^H_INIT^$>0$ or \verb^H_INIT^$<0$, respectively).
The user also needs to specify \verb+DELTA_LAMBDA+ (roughly how far from the initial point to look for the neighboring point) to control this bootstrapping process.

Finally, the user can specify an integer parameter \verb+VERBOSE+ to control output generated at run-time.
No output is generated unless the parameter \verb+VERBOSE+ is positive;
With \verb+VERBOSE>=1+, the master process displays diagnostic messages to standard output as the algorithm progresses.
When \verb+VERBOSE>=2+, the master process also creates data files in a user-specified path that display the structure of the rooted trees.
The data files generated are compatible with the \texttt{dot} language for specifying directed graphs with the \textsc{Graphviz} software for visualization of graphs (see \texttt{www.graphviz.org}).
Such graphs are useful for performance-tuning, \ie, for understanding how the data in \texttt{parameters.txt} affect the use of processors.

Many of the core routines in the \pampac library require traversal of the rooted tree in a depth-first (using recursion) or a breadth-first (using a queue) fashion.
The node and queue data structures for managing the tree are documented in \texttt{pampac.h} in the \texttt{src} subdirectory.
This file also describes a data structure for storing and communicating the parameter options parsed from the user's parameter file.
The \pampac library is designed so that users need not know the details of the implementations of these data structures (nor the routines for allocation/deallocation of memory, management of pointers, etc.).
The user need only specify the width and depth of the underlying tree and the related tunable parameters that control the parallel algorithm.

%% file: numerical-examples.tex
\section{Numerical examples}
\label{sec:num}

We present two examples to test the performance of the parallel algorithm. The first example concerns travelling
waves in a $1+1$ dimensional, nonlinear, partial differential equation (PDE) and is distributed with the \pampac library in the \texttt{examples} directory.
The second concerns time-periodic
solutions to the Navier-Stokes equation on a three-dimensional, periodic domain. In both test cases, the
number of unknowns is in the thousands, but the corrector step is computed differently. In the first
test case, we compute and LU-decompose the dense Jacobian matrix, while in the second case, an inexact
Krylov subspace method is used. An additional difference is that the first test case was implemented in 
{\sf C}, but needs to be linked to external numerical libraries, whereas the second is a self-contained
legacy FORTRAN code. In spite of these differences, a similar set of parameters gives rise to similar
speed-up.

% Results for each problem individually
\input{modified_KS-section} % (Lenn)
\input{fluid-section} % (Lenn) 

%% file: modified_KS-section.tex
\subsection{Travelling waves in a modified Kuramoto-Sivashinsky equation}
\label{ex:modified_ks}

The Kuramoto-Sivashinsky equation was derived independently in different contexts,
most notably by \citeN{Kura} for describing phase dynamics in reaction-diffusion systems
and by \citeN{Siva} for describing flame front dynamics. The similarity of its quadratic 
nonlinearity to that of the Navier-Stokes equation makes it a popular test case for numerical methods
for PDEs. We add a second nonlinear term to obtain
\begin{equation}
u_t+u u_x+u_{xx} +\lambda u_{xxxx}-A \sin(u)=0
\label{modKS}
\end{equation} 
The extra nonlinear term breaks the equivariance 
under Galileo boosts so that we can compute families of travelling waves with a uniquely
determined wave speed. The modified equation is still equivariant under translations and under
the reflection symmetry given by $\mathcal{S}:\ (u,x)\rightarrow (-u,-x)$. 

\begin{figure}
\begin{center}
\includegraphics[width=.9\textwidth]{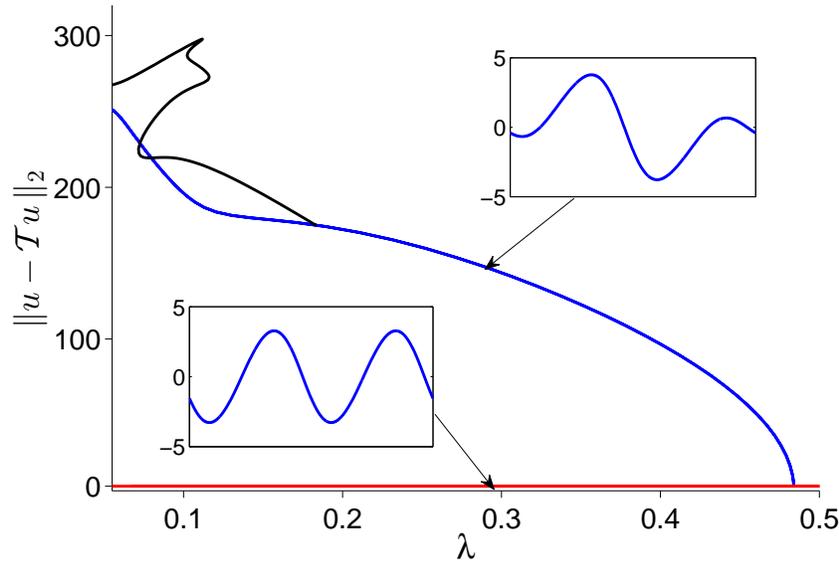}
\end{center}
\caption{Partial bifurcation diagram of the discretization of PDE model \ref{modKS}. In red: branch of equibria 
invariant under symmetries $\mathcal{S}$ and $\mathcal{T}$. In blue: branch of equilibria invariant only
under $\mathcal{S}$. In black the branch of travelling wave solutions shown in detail in Fig.~\ref{KS-curve}. Shown
is a measure for the deviation from invariance under the shift $\mathcal{T}$ versus the viscosity.}
\label{fig:branches}
\end{figure}
\begin{figure}
\begin{center}
\includegraphics[width=.9\textwidth]{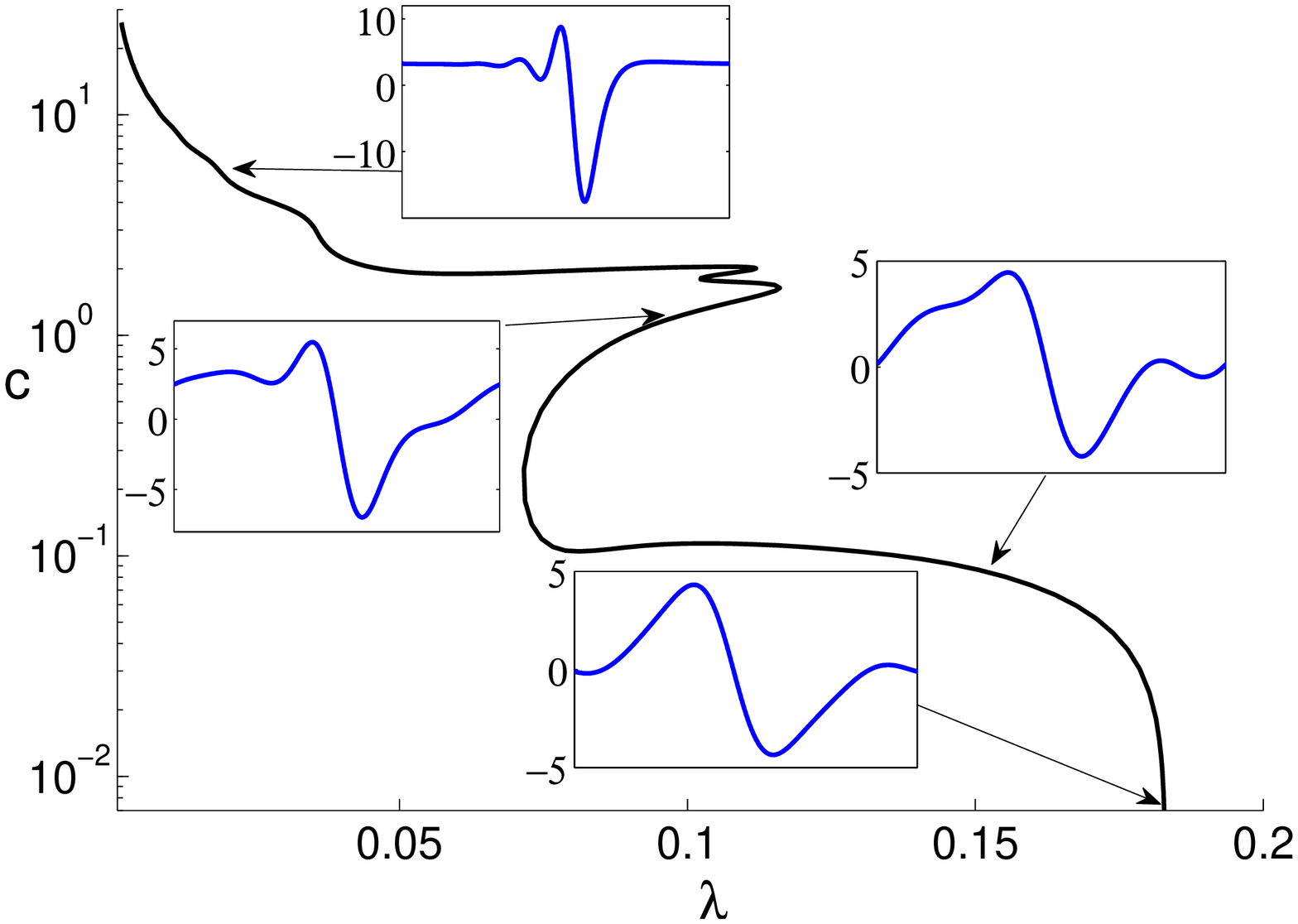}
\end{center}
\caption{The branch of travelling wave solutions to Eq.\ref{modKS} used for testing the parallel continuation
algorithm. Shown is the wave speed, $c$, versus the control parameter, $\lambda$. The inlays show snap shots of the solution
$u(x,t)$ at four points along the continuation curve. At $c=0$, the wave bifurcates
from an equilibrium, which is symmetric under the reflection symmetry. For $\lambda\lesssim 0.01$ the solution becomes
strongly localized.}
\label{KS-curve}
\end{figure}

We consider this PDE on 
a periodic domain, i.e. $u\in C^4(S)$, fixing $A=8.09$ and considering the viscosity, $\lambda>0$, as control parameter.
There is always a trivial solution at $u\equiv 0$, from which equilibria branch off
with increasing wave number for decreasing viscosity. In Fig.~\ref{fig:branches} a branch with wave number two
is shown. This branch is symmetric under the reflection $S$, as well as under the shift $\mathcal{T}$ over half the 
domain. At $\lambda\approx 0.48$, a a family of equilibria branches off in a bifurcation that breaks the translational
symmetry. Subsequently, at $\lambda\approx 0.18$, a family of travelling waves is created, shown in detail in Fig.~\ref{KS-curve}.

This family has a number of fold points and the wave changes its shape rapidly
along the branch. Therefore, its computation by the traditional pseudo-arclength continuation approach suffers
from many failed steps. Moreover, the wave becomes increasingly localized for small values
of the viscosity. In order to resolve it correctly, we need fine discretization,
which will result in time-consuming corrector steps.

We compute the travelling waves as $u(x,t)=w(x-ct)$, which results in the following Boundary Value Problem (BVP):
$$
-c w'+w w'+w''+\lambda w''''-A \sin(w)=0;\ \ \ w(0)=w(2\pi)
$$
The linear terms in this equation are efficiently computed in Fourier space, while the nonlinear terms are
best computed on a regular, periodic grid. Thus, we approximate solutions as
$$
w_j=\sum_{k=0}^{n-1} a_k e^{i k x_j},\ \text{where}\ x_j=\frac{2\pi}{n}j,\ j=0,\,\ldots,\,n-1 
$$
and switch between the representations $\{w_j\}$ and $\{a_k\}$ by the discrete Fourier transform, implemented 
using using the FFTW library \cite{fftw}. The computational cost of each transform is $\mbox{O}(n\ln n)$.
The resolution is fixed to $n=2048$, fine enough to resolve
the solutions shown in Fig.~\ref{KS-curve} and avoid aliasing issues.

The continuation problem now takes the form $F(a_0,\ldots,a_{n-1},c,\lambda)=0$, where $F$ and its derivatives are evaluated
using the pseudo-spectral method outlined above. Since we have an extra unknown, the wave speed $c$, we must add an
extra equation to ensure uniqueness of the corrector steps. We impose that the Newton update steps be orthogonal
to $w_x$, the generator of the symmetry group of translations. This condition ensures that the succesive iterates
under corrector steps do not slide along the $x$-direction.
Because of the second nonlinear term in Eq.~\ref{modKS}, the
Jacobian matrix $\bm{F}_{\bm{x}}$ is dense. The test code takes $\mbox{O}(n^2)$ flops to compute this matrix and 
a further $\mbox{O}(n^3)$ flops to LU-decompose it when solving for the Newton update step.
Thus, the execution of a single Newton step in this test has a complexity of $\mbox{O}(n^3)$ and
takes a few seconds on a single, 2.2GHz CPU, rendering the MPI communication time negligible.
In Table~\ref{constants_mKS} all parameters relevant for the numerical computation are given.
\begin{table}
\tbl{System and algorithm parameters for the numerical experiments with the modified Kuramoto-Sivashinsky equation.
\label{constants_mKS}}{
\begin{tabular}{|c|c|c|}\hline
$A$ & 8.09 & amplitude of second nonlinear term \\
$[\lambda_0,\,\lambda_1]$ & $[0.1828,\,0.001]$ & range of the continuation parameter \\
$n$ & 2048 & \# grid points\\
$\{t_i\}_{i=1}^{3}$ & $\{0.75,\,1,\,2\}$ & step-size multipliers\\
$H_{\rm max}$ & 2000 & maximal step size\\
$H_{\rm min}$ & $10^{-2}$ & minimal step size \\
$H_{\rm init}$ & 100 & initial step size\\
$\nu_{\rm max}$ & 4 & maximal number or Newton iterations\\
$\mu$ & 0.5 & minimal linear redisual decrease\\
$\gamma$ & 2.0 & expected order of residual decrease \\
$r_{\rm max}$ & $5\cdot 10^{-7}$ & tolerance for the nonlinear problem\\\hline
\end{tabular}}
\end{table}

\begin{figure}
\begin{center}
\includegraphics[width=.9\textwidth]{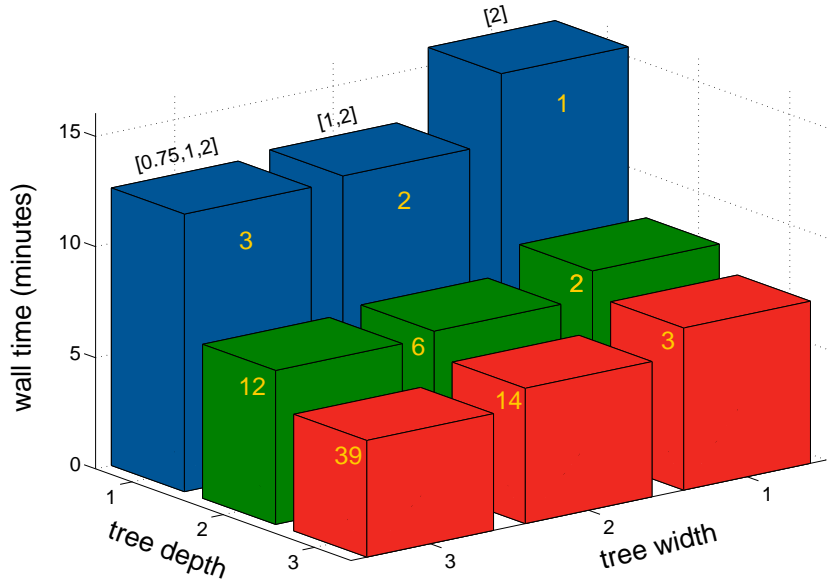}
\end{center}
\caption{Wall time for the computation of the continuation curve shown in Fig.~\ref{KS-curve} using a tree width
and depth up to three. The numbers between brackets denote the multipliers for the step-size, $\{t_i\}_{i=1}^W$ 
and the integer denotes the number of CPUs concurrently processing corrector steps. }
\label{KS-wtime}
\end{figure}

The wall time for the computation of the entire curve shown in Fig.~\ref{KS-curve} using various tree structures
and step-size distributions is shown in Fig.~\ref{KS-wtime}. The maximal speed-up is about a factor of three,
obtained with a tree depth and width of three. In that computation, there is a lot of redundancy as many processors
will be working on very similar approximate solutions. A speed-up of more than a factor of two, however, can be obtained
using as few as three processors with the same step-size multiplication factor. The three depth seems to impact
the efficiency much more strongly than the tree width, but this is likely problem-dependent. Some inital experimentation
will be necesasary for each individual problem to determine a near-optimal strategy for a given number of available processors.

%% file: fluid-section.tex
\subsection{Periodic solutions in a turbulent flow}
\label{ex:turbulence}

The second test problem concerns the continuation of time-periodic solutions to the Navier-Stokes equation for fluid motion.
We consider an incompressible, viscous fluid in a box with periodic boundary conditions in every direction.
In the simulation code, the unknown variables are the Fourier coefficients of the vorticity field, truncated to a finite number. 
Energy is input by keeping fixed in time some coefficients with small wave numbers, corresponding to large spatial scales.
At small spatial scales, energy is dissipated by viscous processes.
The resulting turbulent flow is statistically stationary and exhibits a cascade of energy across spatial scales.
Although the statistical description of this cascade is well developed (see, e.g., \citeN[Ch.~7]{monin}), the dynamics of this process are largely unknown.

In \citeN{veen}, time-periodic solutions of this flow are considered as building blocks of turbulence.
The idea is to study the dynamics and parameter dependence of such building blocks and distill from the results a hypothesis on the dynamical processes that contribute to the energy cascade.
The essential difficulty in computation is that a very large number of degrees of freedom is required to simulate of the flow accurately.
Even for weakly turbulent flow, about $n\simeq10^6$ degrees of freedom are needed; this number increases algebraically with Taylor's microscale
Reynolds number, $\Re_{\lambda}$, which, in turn, increases with decreasing viscosity.
Symmetry arguments reduce the number of degrees of freedom to $n\simeq10^4$, thereby making the study of the resulting \textit{high symmetric flow} feasible \cite{kida}. 

After symmetry reduction, the resulting system is a set of $n$ coupled, nonlinear ODEs with a single parameter, namely the kinematic viscosity, here denoted by $\lambda$, which determines the Reynolds number of the flow.
Time-stepping is done by the pseudo-spectral method, employing the fourth-order Runge-Kutta-Gill scheme.
In the results presented here, the spatial resolution is fixed at $2^7$ points in every direction which, after de-aliasing and symmetry reduction, gives $n=6370$ variables.
The kinematic viscosity is varied in the range $0.0045\geq \lambda \geq 0.0035$, which corresponds to a Taylor microscale Reynolds number $\Re_{\lambda}$ in the range $57 \leq \Re_{\lambda}\leq 68$.
The square  of this dimensionless number can be compared to the commonly used geometric Reynolds number $\Re$.

Periodic solutions are computed as fixed points of an iterated Poincar\'{e} map, i.e., as solutions of a nonlinear system
\begin{equation}
\mathcal{P}^{(k)}(\vec{x},\lambda)-\vec{x}=\vec{F}(\vec{x},\lambda)=\vec{0}.
\label{poincare}
\end{equation}
In \eqref{poincare}, $\vec{x}$ is the vector of Fourier coefficients of the vorticity field and $k$ is the discrete period of the orbit.
The Poincar\'{e} plane of interesection is a coordinate plane on which one of the low wavenumber Fourier coefficients 
equals its time-averaged value.
A solution of discrete period $k=5$ is filtered from turbulent data at the highest viscosity.
At this viscosity, the flow is relatively quiescent.
Subsequently, pseudo-arclength continuation is used to track the periodic solution to the more turbulent regime.
In the correction step of pseudo-arclength continuation, the linear problem associated with the Newton-Raphson iteration is solved in an inner Generalised Minimal Residual (GMRES) iteration \cite{saad}.
The combination of pseudo-arclength continuation with a Krylov subspace iteration is called \textit{Newton-Krylov continuation} and was first implemented by \citeN{sanchez}.
Each linear problem within this continuation takes about twenty inner GMRES iterations to solve.
In turn, each inner GMRES iteration requires integrating a system of $n$ ODEs modeling the flow along an approximately periodic orbit.
On a single, 2.2GHz CPU, one corrector step takes about 13 minutes, again rendering the MPI communication time negligible.
All system and algorithm parameters used to generate the results below are listed in Table~\ref{constants_HS}.
\begin{table}
\tbl{System and algorithm parameters for the numerical experiments with high-symmetric flow.
\label{constants_HS}}{
\begin{tabular}{|c|c|c|}\hline
$[\lambda_0,\,\lambda_1]$ & $[0.0045,\,0.0035]$ & range of the continuation parameter \\
$N$ & $128^3$ & spatial resolution \\
$n$ & 6370 & \# variables\\
$\{t_i\}_{i=1}^{3}$ & $\{0.75,\,1,\,2\}$ & step-size multipliers\\
$h_{\rm i}$ & 0.01 & initial step-size\\
$q$ & 0.5 & minimal redisual decrease is $\|\bm{r}^{(\alpha,\nu_{\alpha})}\|_2< q\|\bm{r}^{(\alpha,\nu_{\alpha}-1)}\|_2$\\
$t_{\rm f}$ & 0.5 & step-size multiplier if all child nodes fail\\
TOL & $10^{-6}$ & tolerance for the nonlinear problem\\
GMRESTOL & $10^{-5}$ & relative residual tolerance for GMRES\\
$\Delta t$ & $5\times 10^{-3}$ & step-size for $4^{\rm th}$ order Runge-Kutta-Gill time stepping\\\hline
\end{tabular}}
\end{table}
\begin{figure}
\begin{center}
\includegraphics[width=.6\textwidth]{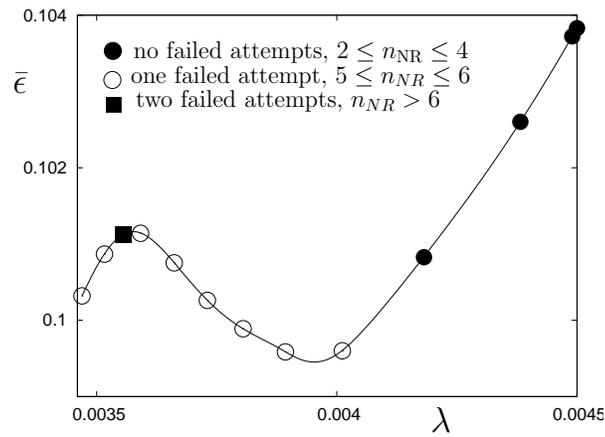}
\caption{Continuation of a periodic solution in high-symmetric flow by serial pseudo-arclength continutation.
Dots, circles, and squares denote points computed after zero, one and two failed correction steps, respectively.}
\label{serial_high_symm}
\end{center}
\end{figure}

\begin{figure}
\begin{center}
\includegraphics[width=.9\textwidth]{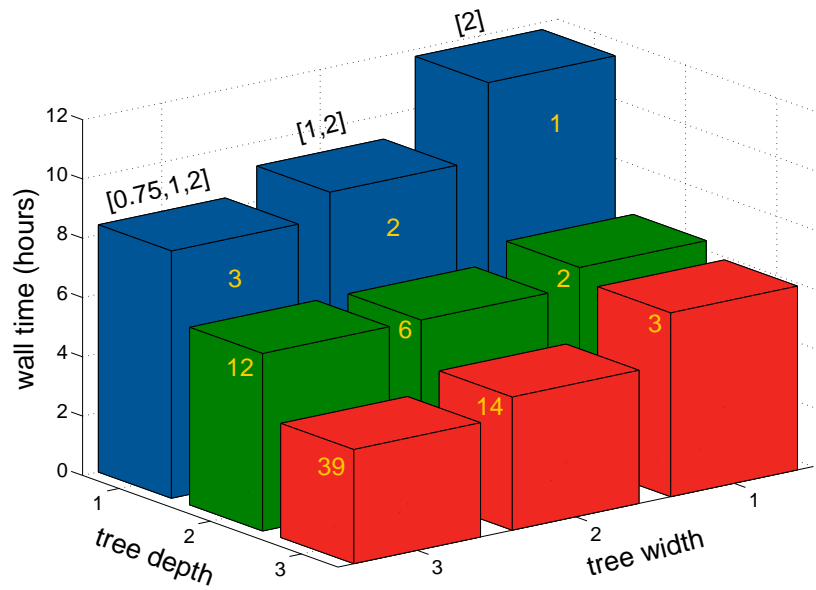}
\caption{Wall time in hours versus tree width and depth. The numbers between brackets are the step length multiplication factors $\{t_i\}_{i=1}^W$ and the integer
on each data bar is the number of CPUs concurrently processing corrector steps.}
\label{wtime_fluid}
\end{center}
\end{figure}

Figure~\ref{serial_high_symm} shows the numerical curve produced by pseudo-arclength continuation.
The points shown are computed using a naive scheme for step-size control: the step-size is doubled after each successful step, halved after each failed step.
Twelve points are computed along the curve at the cost of fifty-five Newton-Krylov iterations.
Out of these fifty-five correction steps, fifteen steps were rejected; thus, over a quarter of the time (as measured on a wall-clock) was wasted.
A more careful strategy based on a local estimate of the curvature yields eighteen points computed on the numerical curve at a cost of forty-nine Newton iterations with four rejected steps.
In Fig.~\ref{wtime_fluid} the wall time required by \pampac is shown for tree configurations up to width and depth three. 
As compared to the first test case, the tree width, i.e. the number of different step lengths attempted in parallel,
is of greater influence. However, a speed-up by a factor of two is again obtained using three CPUs and the maximal speed-up
by a factor of three is obtained on 39 CPUs.

%% file: conclusion.tex
\section{conclusion}\label{sec:conclusion}

In the computation of parametrised solutions to discretised PDEs, most elements have been optimised
for efficiency. For instance, when studying Navier--Stokes flow, we often use pseudo-spectral time-stepping
in combination with fast Fourier transforms. The linear systems we must solve to find Newton update steps
are handled by Krylov subspace methods, whose convergence can be sped up by a host of preconditioners.
In contrast, the continuation algorithm, which forms the outer loop of the computation, is essentially
the same as that used for small sets of ODEs. For such small systems, selecting an unnecessarily small
step-size, or a overly large step-size that leads to diverging corrector steps, will cost seconds
or minutes of computation time. For large systems, this may cost days or weeks.

In this paper, we have presented an elegant, recursive algorithm that combines two strategies for
aggressively optimising the step-size and minimising the computation time. The first is to try 
several step-sizes in parallel, and the second is to predict new solutions from a sequence of
corrector steps, before this sequence has converged. 

Two test cases, different in the type of
solution computed, the linear solving and the implementation, demonstrate that the continuation 
can be sped up by a factor of two using only three processors, and by a factor of three using 
thirty nine. Since multi-core processors are a standard feature of new desktop computers, and 
cluster computers are available in many places, we expect that our algorithm will be useful
to many researchers working on nonlinear problems with many unknowns.